\UseRawInputEncoding
\documentclass[oneside, 11pt]{amsart}
\usepackage{amsmath}
\usepackage{amssymb}
\usepackage{dsfont}
\usepackage{color}
\usepackage[usenames,dvipsnames,x11names,svgnames]{xcolor}
\usepackage[colorlinks=true,linkcolor=black,citecolor=black, urlcolor=Blue]{hyperref}
\usepackage{amsthm}
\usepackage{amscd}
\usepackage{geometry}
\usepackage{mathrsfs}
\usepackage{mathtools}

\newtheorem{theorem}{Theorem}

\newtheorem{lemma}{Lemma}
\newtheorem{corollary}{Corollary}

\newtheorem{problem}{Problem}

\theoremstyle{definition}
\newtheorem{definition}{Definition}
\newtheorem*{remark}{Remark}

\newcommand{\be}{\begin{equation*}}
\newcommand{\ee}{\end{equation*} }

\newcommand{\ben}{\begin{equation}}
\newcommand{\een}{\end{equation} }

\newcommand{\bs}{\begin{split}}
\newcommand{\es}{\end{split}}

\newcommand{\bmu}{\begin{multline*}}
\newcommand{\emu}{\end{multline*}}

\newcommand{\bmun}{\begin{multline}}
\newcommand{\emun}{\end{multline}}

\begin{document}



\title[]{On all real zeros for a new class of the even entire function}

\author{Xiao-Jun Yang$^{1,2}$}

\email{dyangxiaojun@163.com; xjyang@cumt.edu.cn}

\address{$^{1}$ School of Mathematics,and State Key Laboratory for Geo-Mechanics and
Deep Underground Engineering, China University of Mining and Technology, Xuzhou 221116, China}
\address{$^{2}$ Department of Mathematics, Faculty of Science, King Abdulaziz University P.O. Box 80257, Jeddah 21589, Saudi Arabia}

\begin{abstract}
In this article we propose a new class of the even entire function connected
with the product and series with the real coefficients. We address a
sufficient condition for all real zeros for it. As a typical example, we
give an answer to the problem of Lagarias and Montague. We suggest the open
problems for the class of the even entire function.
\end{abstract}

\maketitle
\tableofcontents
%%%%%%%%%%%%%%%%%%%%%%%%%%%%%%%%%%%%%%%%%%%%%%%%%%%%%%%%%%%%%%%%%%%%%%%%%%%%%%%%%%%%%%%%%%%
%\\maketitle                                %！！Show title
%\tableofcontents                             %！！Content(1)
     %\\section{}                             %！！1 title
      %\subsection{}                         %！！2 title
       %\subsubsection{}                     %！！3 title
%%%%%%%%%%%%%%%%%%%%%%%%%%%%%%%%%%%%%%%%%%%%%%%%%%%%%%%%%%%%%%%%%%%%%%%%%%%%%%%%%%%%%%%%%%%%%%%%%%%%%%

%%%%%%%%%%%%%%%%%%%%%%%%%%%%%%%%%%%%%%%%%%%%%%%%%%%%%%%%%%%%%%%%%%%%%%%%%%%%%%%%%%%%%%%%%%%%%%%%%%%%%%
%%%%%%%%%%%%%%%%%%%%%%%%%%%%%%%%%%%%%%%%%%%%%%%%%%%%%%%%%%%%%%%%%%%%%%%%%%%%%%%%%%%%%%%%%%%%%%%%%%%%%%
\section{Introduction} {\label{sec:1}}

The theory of the entire functions has played an important role in the study
of the behaviors of the zeros of the functions of the real and complex
variables (see for instance \cite{1,2} and the references therein). The entire
functions have the popular applications of the polynomials of the entire
functions with the distributions of their zeros in the fields of
mathematics, physics, and engineering (see \cite{3,4,5,6}).
Let $\mathbb{C}$ and $\mathbb{R}$ be the sets of the complex and real numbers,
respectively, and let $i=\sqrt {-1} $. In 2015 Csordas the entire function
$\Xi \left( \tau \right)$ has a nice connection with the Fourier transform \cite{7}
\begin{equation}
\label{eq1}
\begin{array}{l}
\Omega \left( \tau \right)=\int\limits_{-\infty }^\infty {\hbar \left( \ell
\right)e^{i\tau\ell }d\ell } =\int\limits_0^\infty {\hbar \left( \ell
\right)\cos \left( {\tau\ell } \right)d\ell } ,
\end{array}
\end{equation}
where $\hbar \left( \ell \right)$ is a positive-value kernel in the domain
$\ell \in \mathbb{R}$ and $t\in \mathbb{C}$. In fact, the integral equation (\ref{eq1}) was
proposed by Jensen in 1913 \cite{8} and developed by P\'{o}lya in 1927 \cite{9,10} due
to the connection with the Riemann zeta function \cite{11}. In Bruijn's 1950
paper, he considered (\ref{eq1}) to set up the zeroes of the Ramanujan zeta
function \cite{12}.
Jensen \cite{8} proposed that Riemann xi function $\Xi \left( t \right)$ has the
Fourier cosine integral representation as follows:
\begin{equation}
\label{eq2}
\begin{array}{l}
\Xi \left( \tau \right)=2\int\limits_0^\infty {g\left( \ell \right)\cos \left(
{\tau \ell } \right)d\ell } ,
\end{array}
\end{equation}
where
\begin{equation}
\label{eq3}
\begin{array}{l}
g\left( \ell \right)=\sum\limits_{n=1}^\infty {\left( {4\pi ^2 n
^4e^{\frac{9}{2}\ell }-6\pi n ^2e^{\frac{5}{2}\ell }} \right)} e^{-\pi
n^2e^{2\ell }}.
\end{array}
\end{equation}
In 2019 Griffin, Ono, Rolen and Zagier \cite{13} made a big contribution for the real zeros of
(\ref{eq2}). In 1976 Newman \cite{14} introduced the a family of the Fourier
cosine integral
\begin{equation}
\label{eq4}
\begin{array}{l}
H \left( \tau \right)=2\int\limits_0^\infty {e^{-\nu \ell ^2}g\left( \ell
\right)\cos \left( {\tau\ell } \right)d\ell } ,
\end{array}
\end{equation}
where $\nu \in \mathbb{R}$ is the de Bruijn-Newman constant.
In 2020 Rodgers and Tao \cite{15} reported that $\nu $ is non-negative can be considered as
the non-negative value. Ki, Kim and Lee \cite{16} suggested that the de
Bruijn-Newman constant $\nu <1/2$.
In 2011 Lagarias and Montague \cite{17} suggested the Fourier sine integral
\begin{equation}
\label{eq5}
\begin{array}{l}
\Phi \left( \tau \right)=2\int\limits_0^\infty {g\left( \ell \right)\ell
^{-1}\sin \left( {\tau\ell } \right)d\ell }.
\end{array}
\end{equation}
The problem of Lagarias and Montague is a conjecture given by Lagarias and Montague \cite{17},
which states all zeros of the
function $\Phi \left( \tau \right)$ are real. As a generalization of the work in
\cite{18} to solve the problem of Jensen \cite{8}, the main target of the present paper is to consider a family of the
even entire function as follows:

\begin{definition}
An even entire function $\aleph \left( \tau \right)$ of order $\gamma =1$
for $\tau \in \mathbb{C}$, represented by the series
\begin{equation}
\label{eq6}
\begin{array}{l}
\aleph \left( \tau \right)=\sum\limits_{m=0}^\infty {\alpha _m \tau ^{2m}},
\end{array}
\end{equation}
is said to be in the class $\wp $, written as $\aleph \in \wp $, if $\aleph
\left( \tau \right)$ can be expressed as
\begin{equation}
\label{eq7}
\begin{array}{l}
\aleph \left( \tau \right)=\aleph \left( 0 \right)\prod\limits_{\vartheta _k
} {\left( {1-\frac{\tau i}{\vartheta _k }} \right)} ,
\end{array}
\end{equation}
where $\vartheta _k \ne 0$ run the zeros of $\aleph \left( \tau \right)$,
$\alpha _m =\left( {-1} \right)^mB_m $ are the coefficients for $\aleph
\left( \tau \right)$ with $B_m >0$, the series
\begin{equation}
\label{eq8}
\begin{array}{l}
p=\sum\limits_{k=1}^\infty {\frac{1}{|\vartheta _k|^2 }},
\end{array}
\end{equation}
is convergent, and
\begin{equation}
\label{eq9}
\begin{array}{l}
\aleph \left( 0 \right)\ne 0.
\end{array}
\end{equation}
\end{definition}
It is obvious that the even entire function $\aleph \left( \tau \right)$ has
the infinite zeros because it is not constant and finite polynomials.
The outline of this paper is given as follows.
In Section \ref{sec:2} we introduce the equivalent idea for the class of the even
entire function. In Section \ref{sec:3} we prove:
\begin{theorem}
\label{Th1}
Let $\aleph \in \wp $. Then all of the zeros of $\aleph \left( \tau
\right)$ are real.
\end{theorem}
In Section \ref{sec:4} we give a detailed answer for the problem of Lagarias and
Montague.

\section{A special class of the even entire functions} {\label{sec:2}}

\subsection{Change the product of $\aleph \left( \tau \right)$}
We start with the following result:

\begin{lemma}
\label{Le1}
If $\aleph \in \wp $, then there is
\begin{equation}
\label{eq10}
\begin{array}{l}
\aleph \left( \tau \right)=\aleph \left( 0 \right)\prod\limits_{\Im \left(
{\vartheta _k } \right)>0} {\left( {1+\frac{\tau ^2}{\vartheta _k^2 }}
\right)} .
\end{array}
\end{equation}
\end{lemma}

\begin{proof}
Because $\aleph \left( \tau \right)$ is an even function, we present
\begin{equation}
\label{eq11}
\begin{array}{l}
\aleph \left( {\vartheta _k } \right)=\aleph \left( {-\vartheta _k }
\right).
\end{array}
\end{equation}
By (\ref{eq11}), the function (\ref{eq7}) can be written as
\begin{equation}
\label{eq12}
\begin{array}{l}
 \aleph \left( \tau \right)
 =\aleph \left( 0 \right)\prod\limits_{\vartheta
_k } {\left( {1-\frac{\tau i}{\vartheta _k }} \right)}
 =\aleph \left( 0 \right)\prod\limits_{\Im \left( {\vartheta _k } \right)>0}
{\left( {1-\frac{\tau i}{\vartheta _k }} \right)\left( {1-\frac{\tau
i}{\left( {-\vartheta _k } \right)}} \right)}
 =\aleph \left( 0 \right)\prod\limits_{\Im \left( {\vartheta _k } \right)>0}
{\left( {1+\frac{\tau ^2}{\vartheta _k^2 }} \right)} .
 \end{array}
\end{equation}
We thus finish the proof.
\end{proof}

\subsection{Structure the function $\psi \left( \vartheta \right)$}
Now, let us give the special class of the even entire functions.

We first set
\begin{equation}
\label{eq13}
\begin{array}{l}
x=-i\vartheta ,
\end{array}
\end{equation}
where $\vartheta \in \mathbb{C}$.

Substituting (\ref{eq13}) back into (\ref{eq7}), we obtain
\begin{equation}
\label{eq14}
\begin{array}{l}
\aleph \left( {-i\vartheta } \right)=\sum\limits_{m=0}^\infty {B_m x^{2m}}
\end{array}
\end{equation}
and
\begin{equation}
\label{eq15}
\begin{array}{l}
\aleph \left( {-i\vartheta } \right)=\aleph \left( 0
\right)\prod\limits_{\vartheta _k } {\left( {1-\frac{\vartheta }{\vartheta
_k }} \right)} .
\end{array}
\end{equation}
Let us define the function $\psi \left( \vartheta \right)$ by
\begin{equation}
\label{eq16}
\begin{array}{l}
\psi \left( \vartheta \right):=\aleph \left( {-i\vartheta } \right).
\end{array}
\end{equation}
Then (\ref{eq14}) and (\ref{eq15}) can be rewritten as
\begin{equation}
\label{eq17}
\begin{array}{l}
\psi \left( \vartheta \right)=\sum\limits_{m=0}^\infty {B_m \vartheta ^{2m}}
\end{array}
\end{equation}
and
\begin{equation}
\label{eq18}
\begin{array}{l}
\psi \left( \vartheta \right)=\psi \left( 0 \right)\prod\limits_{\vartheta
_k } {\left( {1-\frac{\vartheta }{\vartheta _k }} \right)} ,
\end{array}
\end{equation}
respectively.

\begin{lemma}
\label{Le2}
There exists the functional equation
\begin{equation}
\label{eq19}
\begin{array}{l}
\psi \left( \vartheta \right)=\psi \left( {-\vartheta } \right).
\end{array}
\end{equation}
\end{lemma}

\begin{proof}
From (\ref{eq17}) we have
\begin{equation}
\label{eq20}
\begin{array}{l}
\psi \left( {-\vartheta } \right)=\sum\limits_{m=0}^\infty {B_m \left(
{-\vartheta } \right)^{2m}} =\sum\limits_{m=0}^\infty {B_m \vartheta ^{2m}}.
\end{array}
\end{equation}
Thus, this is the required result.
\end{proof}

\begin{lemma}
\label{Le3}
There exists the identity
\begin{equation}
\label{eq21}
\begin{array}{l}
\psi \left( 0 \right)\prod\limits_{\vartheta _k } {\left( {1-\frac{\vartheta
}{\vartheta _k }} \right)} =\psi \left( 0 \right)\prod\limits_{\Im \left(
{\vartheta _k } \right)>0} {\left( {1-\frac{\vartheta ^2}{\vartheta _k^2 }}
\right)} .
\end{array}
\end{equation}
\end{lemma}

\begin{proof}
By using Lemma \ref{Le2}, we get
\begin{equation}
\label{eq22}
\begin{array}{l}
 \psi \left( \vartheta \right)=\psi \left( 0 \right)\prod\limits_{\vartheta
_k } {\left( {1-\frac{\vartheta }{\vartheta _k }} \right)}
 =\psi \left( 0 \right)\prod\limits_{\Im \left( {\vartheta _k } \right)>0}
{\left( {1-\frac{\vartheta }{\vartheta _k }} \right)} \left(
{1+\frac{\vartheta }{\vartheta _k }} \right)
 =\psi \left( 0 \right)\prod\limits_{\Im \left( {\vartheta _k } \right)>0}
{\left( {1-\frac{\vartheta ^2}{\vartheta _k^2 }} \right)} .
 \end{array}
\end{equation}
We thus complete the proof.
\end{proof}

Combining (\ref{eq17}), (\ref{eq18}) and Lemma \ref{Le3}, we obtain
\begin{equation}
\label{eq23}
\begin{array}{l}
 \psi \left( \vartheta \right)=\sum\limits_{m=0}^\infty {B_m \vartheta
^{2m}}
 =\psi \left( 0 \right)\prod\limits_{\vartheta _k } {\left(
{1-\frac{\vartheta }{\vartheta _k }} \right)}
 =\psi \left( 0 \right)\prod\limits_{\Im \left( {\vartheta _k } \right)>0}
{\left( {1-\frac{\vartheta ^2}{\vartheta _k^2 }} \right)} .
 \end{array}
\end{equation}
Considering the given condition that $\aleph \left( \tau \right)$ is an even function of order
$\gamma =1$, the argument in Levin's book
(see Theorem 13 in \cite{1}, p.24) said that there exists $\varepsilon >0$ such that
\begin{equation}
\label{eq24}
\begin{array}{l}
\sum\limits_{k=1}^\infty {\frac{1}{\left| {-i\vartheta _k }
\right|^{1+\varepsilon }}} =\sum\limits_{k=1}^\infty {\frac{1}{\left|
{\vartheta _k } \right|^{1+\varepsilon }}}
\end{array}
\end{equation}
is convergent.
Taking $\varepsilon =1$ in (\ref{eq24}), we know that
\begin{equation}
\label{eq25}
\begin{array}{l}
\sum\limits_{k=1}^\infty {\frac{1}{\left| {\vartheta _k } \right|^2}}
\end{array}
\end{equation}
is convergent.

\section{The proof of Theorem \ref{Th1}} {\label{sec:3}}

We now begin with the proof of the Theorem \ref{Th1}. We need to
divide it into six steps.

\subsection{Set up a class of $\psi \left( \vartheta \right)$}

By (\ref{eq6}) and (\ref{eq7}), we structure
\begin{equation}
\label{eq26}
\begin{array}{l}
\sum\limits_{m=0}^\infty {B_m \vartheta ^{2m}} =\psi \left( 0
\right)\prod\limits_{\vartheta _k } {\left( {1-\frac{\vartheta }{\vartheta
_k }} \right)} .
\end{array}
\end{equation}
Owing to Lemma \ref{Le3}, we easily demonstrate that
\begin{equation}
\label{eq27}
\begin{array}{l}
\psi \left( \vartheta \right)=\psi \left( 0 \right)\prod\limits_{\vartheta
_k } {\left( {1-\frac{\vartheta }{\vartheta _k }} \right)} =\psi \left( 0
\right)\prod\limits_{\Im \left( {\vartheta _k } \right)>0} {\left(
{1-\frac{\vartheta ^2}{\vartheta _k^2 }} \right)} .
\end{array}
\end{equation}
Combining (\ref{eq26}) into (\ref{eq27}), we show that
\begin{equation}
\label{eq28}
\begin{array}{l}
\psi \left( \vartheta \right)=\sum\limits_{m=0}^\infty {B_m \vartheta ^{2m}}
=\psi \left( 0 \right)\prod\limits_{\Im \left( {\vartheta _k } \right)>0}
{\left( {1-\frac{\vartheta ^2}{\vartheta _k^2 }} \right)} .
\end{array}
\end{equation}

\subsection{Find a class of $\overline {\psi \left( \vartheta \right)} $}

Let $\overline {\psi \left( \vartheta \right)} $, $\overline \vartheta $ and
$\overline \vartheta _k $ be the complex conjugates of $\psi \left(
\vartheta \right)$, $\vartheta $ and $\vartheta _k $, respectively.

By (\ref{eq28}), the function $\overline {\psi \left( \vartheta \right)} $ can be expressed as
\begin{equation}
\label{eq29}
\begin{array}{l}
\overline {\psi \left( \vartheta \right)} =\overline {\left[
{\sum\limits_{m=0}^\infty {B_m \vartheta ^{2m}} } \right]}
=\sum\limits_{m=0}^\infty {B_m \overline \vartheta ^{2m}}
\end{array}
\end{equation}
due to $B_m >0$.

It follows from (\ref{eq29}) that
\begin{equation}
\label{eq30}
\begin{array}{l}
\overline {\psi \left( \vartheta \right)} =\psi \left( {\overline \vartheta
} \right).
\end{array}
\end{equation}
Making use of (\ref{eq28}), we have
\[
\begin{array}{l}
\psi \left( \vartheta \right)=\psi \left( 0 \right)\prod\limits_{\Im \left(
{\vartheta _k } \right)>0} {\left( {1-\frac{\vartheta ^2}{\vartheta _k^2 }}
\right)},
\end{array}
\]
which implies from (\ref{eq30}) that
\begin{equation}
\label{eq31}
\begin{array}{l}
\overline {\psi \left( \vartheta \right)} =\psi \left( 0
\right)\prod\limits_{\Im \left( {\vartheta _k } \right)>0} {\left(
{1-\frac{\overline \vartheta ^2}{\vartheta _k^2 }} \right)}.
\end{array}
\end{equation}
Similarly, by (\ref{eq28}), we have
\[
\begin{array}{l}
\psi \left( \vartheta \right)=\psi \left( 0 \right)\prod\limits_{\Im \left(
{\vartheta _k } \right)>0} {\left( {1-\frac{\vartheta ^2}{\vartheta _k^2 }}
\right)}
\end{array}
\]
such that
\begin{equation}
\label{eq32}
\begin{array}{l}
 \overline {\psi \left( \vartheta \right)}
 =\overline {\left[ {\psi \left( 0
\right)\prod\limits_{\Im \left( {\vartheta _k } \right)>0} {\left(
{1-\frac{\vartheta ^2}{\vartheta _k^2 }} \right)} } \right]}
=\psi \left( 0
\right)\overline {\left[ {\prod\limits_{\Im \left( {\vartheta _k }
\right)>0} {\left( {1-\frac{\vartheta ^2}{\vartheta _k^2 }} \right)} }
\right]}
 =\psi \left( 0 \right)\prod\limits_{\Im \left( {\vartheta _k } \right)>0}
{\left( {1-\frac{\overline \vartheta ^2}{\overline \vartheta _k^2 }}
\right)}
 \end{array}
\end{equation}
since
\begin{equation}
\label{eq33}
\begin{array}{l}
\psi \left( 0 \right)=B_0 >0,
\end{array}
\end{equation}
proved that one takes $\vartheta =0$ into (\ref{eq28}).

Combining (\ref{eq31}) and (\ref{eq32}), we suggest
\begin{equation}
\label{eq34}
\begin{array}{l}
\psi \left( 0 \right)\prod\limits_{\Im \left( {\vartheta _k } \right)>0}
{\left( {1-\frac{\overline \vartheta ^2}{\vartheta _k^2 }} \right)} =\psi
\left( 0 \right)\prod\limits_{\Im \left( {\vartheta _k } \right)>0} {\left(
{1-\frac{\overline \vartheta ^2}{\overline \vartheta _k^2 }} \right)} .
\end{array}
\end{equation}
In view of (\ref{eq30}) and (\ref{eq34}), we demonstrate the identity
\begin{equation}
\label{eq35}
\begin{array}{l}
 \overline {\psi \left( \vartheta \right)} =\psi \left( {\overline \vartheta
} \right)
 =\psi \left( 0 \right)\prod\limits_{\Im \left( {\vartheta _k } \right)>0}
{\left( {1-\frac{\overline \vartheta ^2}{\vartheta _k^2 }} \right)}
 =\psi \left( 0 \right)\prod\limits_{\Im \left( {\vartheta _k } \right)>0}
{\left( {1-\frac{\overline \vartheta ^2}{\overline \vartheta _k^2 }}
\right)} .
 \end{array}
\end{equation}

\subsection{Propose a class of $\psi \left( t \right)$}

Putting
\begin{equation}
\label{eq36}
\begin{array}{l}
\vartheta =t\in \mathbb{R}
\end{array}
\end{equation}
into (\ref{eq29}), we find that
\begin{equation}
\label{eq37}
\begin{array}{l}
\overline {\psi \left( t \right)} =\overline {\left[
{\sum\limits_{m=0}^\infty {B_m t^{2m}} } \right]} =\sum\limits_{m=0}^\infty
{B_m \overline t ^{2m}} =\sum\limits_{m=0}^\infty {B_m t^{2m}} ,
\end{array}
\end{equation}
which implies from (\ref{eq30}) that
\begin{equation}
\label{eq38}
\begin{array}{l}
\overline {\psi \left( t \right)} =\psi \left( {\overline t } \right)=\psi
\left( t \right).
\end{array}
\end{equation}
With (\ref{eq31}) and (\ref{eq38}), we arrive at
\begin{equation}
\label{eq39}
\begin{array}{l}
\overline {\psi \left( t \right)} =\psi \left( 0 \right)\prod\limits_{\Im
\left( {\vartheta _k } \right)>0} {\left( {1-\frac{\overline t ^2}{\vartheta
_k^2 }} \right)} =\psi \left( 0 \right)\prod\limits_{\Im \left( {\vartheta
_k } \right)>0} {\left( {1-\frac{t^2}{\vartheta _k^2 }} \right)} .
\end{array}
\end{equation}
With the aid of (\ref{eq39}), we may get
\begin{equation}
\label{eq40}
\begin{array}{l}
\psi \left( t \right)=\psi \left( 0 \right)\prod\limits_{\Im \left(
{\vartheta _k } \right)>0} {\left( {1-\frac{t^2}{\vartheta _k^2 }} \right)}.
\end{array}
\end{equation}
In a similar way, by using (\ref{eq32}) and (\ref{eq38}), we show
\begin{equation}
\label{eq41}
\begin{array}{l}
 \overline {\psi \left( t \right)}
 =\psi \left( 0 \right)\prod\limits_{\Im
\left( {\vartheta _k } \right)>0} {\left( {1-\frac{\overline t ^2}{\overline
\vartheta _k^2 }} \right)}
 =\psi \left( 0 \right)\prod\limits_{\Im \left( {\vartheta _k } \right)>0}
{\left( {1-\frac{t^2}{\overline \vartheta _k^2 }} \right)} .
 \end{array}
\end{equation}
With (\ref{eq38}) and (\ref{eq41}), we obtain
\begin{equation}
\label{eq42}
\begin{array}{l}
\psi \left( t \right)=\psi \left( 0 \right)\prod\limits_{\Im \left(
{\vartheta _k } \right)>0} {\left( {1-\frac{t^2}{\overline \vartheta _k^2 }}
\right)} .
\end{array}
\end{equation}
Combining (\ref{eq40}) and (\ref{eq42}), we suggest
\begin{equation}
\label{eq43}
\begin{array}{l}
 \psi \left( t \right)
 =\sum\limits_{m=0}^\infty {B_m t^{2m}}
 =\psi \left( 0
\right)\prod\limits_{\Im \left( {\vartheta _k } \right)>0} {\left(
{1-\frac{t^2}{\vartheta _k^2 }} \right)}
 =\psi \left( 0 \right)\prod\limits_{\Im \left( {\vartheta _k } \right)>0}
{\left( {1-\frac{t^2}{\overline \vartheta _k^2 }} \right)} .
 \end{array}
\end{equation}

\subsection{Consider the convergence of $\psi \left( \vartheta \right)$}

Considering the fact $B_m >0$ and using (\ref{eq43}), we have
\begin{equation}
\label{eq44}
\psi \left( 1 \right)=\sum\limits_{m=0}^\infty {B_m}
\end{equation}
such that
\begin{equation}
\label{eq45}
\psi \left( 1\right)>0.
\end{equation}

With use of (\ref{eq45}) and (\ref{eq43}), we find that
\begin{equation}
\label{eq46}
\begin{array}{l}
\psi \left( 1 \right)=\psi \left( 0 \right)\prod\limits_{\Im \left( {\vartheta _k } \right)>0}
{\left( {1-\frac{1}{\vartheta _k^2 }} \right)}>0,
\end{array}
\end{equation}

\begin{equation}
\label{eq47}
\begin{array}{l}
\psi \left( 1 \right)=\psi \left( 0 \right)\prod\limits_{\Im \left( {\vartheta _k } \right)>0}
{\left( {1-\frac{1}{\overline \vartheta _k^2 }} \right)}>0,
\end{array}
\end{equation}
and
\begin{equation}
\label{eq48}
\begin{array}{l}
\psi \left( 0 \right)\prod\limits_{\Im \left( {\vartheta _k } \right)>0}
{\left( {1-\frac{1}{\vartheta _k^2 }} \right)} =\psi \left( 0
\right)\prod\limits_{\Im \left( {\vartheta _k } \right)>0} {\left(
{1-\frac{1}{\overline \vartheta _k^2 }} \right)} .
\end{array}
\end{equation}
It follows from (\ref{eq24}) that the products of (\ref{eq48}) are absolutely convergent.
This implies that the products of (\ref{eq48}) are convergent.

Titchmarsh (see \cite{19}, p.14-15) argued that (\ref{eq48}) is convergent
if and only if
\begin{equation}
\label{eq49}
\sum\limits_{k=1}^\infty {\frac{1}{\vartheta _k^2 }}
\end{equation}
and
\begin{equation}
\label{eq50}
\sum\limits_{k=1}^\infty {\frac{1}{\overline \vartheta _k^2 }}
\end{equation}
are convergent and there always exists
\begin{equation}
\label{eq51}
\sum\limits_{k=1}^\infty {\frac{1}{\overline \vartheta _k^2 }}
=\sum\limits_{k=1}^\infty {\frac{1}{\vartheta _k^2 }}.
\end{equation}
Applying (\ref{eq49}) and (\ref{eq50}), we see that (\ref{eq28}) and (\ref{eq43}) are convergent.

\subsection{Present the identity $\Re \left( {\vartheta _k } \right)=0$}

Since (\ref{eq49}) and (\ref{eq50}) are valid, (\ref{eq51}) can be
written as
\begin{equation}
\label{eq54}
\begin{array}{l}
 \sum\limits_{k=1}^\infty {\frac{1}{\vartheta _k^2 }}
-\sum\limits_{k=1}^\infty {\frac{1}{\overline \vartheta _k^2 }}
=0.
 \end{array}
\end{equation}
From (\ref{eq54}) we obtain
\begin{equation}
\label{eq55}
\overline \vartheta _k^2 -\vartheta _k^2 =0
\end{equation}
or, alternatively,
\begin{equation}
\label{eq56}
\begin{array}{l}
\left( {\overline \vartheta _k -\vartheta _k } \right)\left( {\vartheta _k
+\overline \vartheta _k } \right)=0.
\end{array}
\end{equation}
Since $\Im \left( {\vartheta _k } \right)>0$ and $\vartheta _k -\overline
\vartheta _k =2i\Im \left( {\vartheta _k } \right)$, we have from (\ref{eq56})
that
\begin{equation}
\label{eq57}
\begin{array}{l}
 \vartheta _k +\overline \vartheta _k =2\Re \left( {\vartheta _k } \right) =0,
 \end{array}
\end{equation}
or, alternatively,
\begin{equation}
\label{eq58}
\begin{array}{l}
 \Re \left( {\vartheta _k } \right) =0.
 \end{array}
\end{equation}
Taking
\begin{equation}
\label{eq59}
\begin{array}{l}
 \left| {\Im \left( {\vartheta _k } \right)} \right|=\sigma _k >0
 \end{array}
\end{equation}
and substituting (\ref{eq58}) back into (\ref{eq43}), we obtain
\begin{equation}
\label{eq60}
\begin{array}{l}
 \psi \left( t \right)
 =\psi \left( 0 \right)\prod\limits_{\Im \left(
{\vartheta _k } \right)>0} {\left( {1-\frac{t^2}{\vartheta _k^2 }} \right)}
=\psi \left( 0 \right)\prod\limits_{\sigma _k } {\left( {1+\frac{\vartheta
^2}{\sigma _k^2 }} \right)}
 =\psi \left( 0 \right)\prod\limits_{k=1}^\infty {\left( {1+\frac{\vartheta
^2}{\sigma _k^2 }} \right)}.
 \end{array}
\end{equation}

Adopting (\ref{eq55}), (\ref{eq58}) and (\ref{eq59}), we have
\begin{equation}
\label{eq61}
\begin{array}{l}
 \vartheta _k^2 =\overline \vartheta _k^2
 ={-\sigma _k^2}.
 \end{array}
\end{equation}

It follows from (\ref{eq51}) that
\begin{equation}
\label{eq63}
\begin{array}{l}
 p=\sum\limits_{k=1}^\infty {\frac{1}{\vartheta _k^2 }}
 =\sum\limits_{k=1}^\infty {\frac{1}{\overline \vartheta _k^2 }}
 =-\sum\limits_{k=1}^\infty {\frac{1}{\sigma _k^2 }}.
 \end{array}
\end{equation}

From (\ref{eq25}) we have

\begin{equation}
\label{eq3.401}
\begin{array}{l}
 \sum\limits_{k=1}^\infty {\frac{1}{\left| {\vartheta _k } \right|^2}}
 =\sum\limits_{k=1}^\infty {\frac{1}{\sigma _k^2 }}=-p.
 \end{array}
\end{equation}

Consequently, (\ref{eq58}) is true.

Combining (\ref{eq23}) and (\ref{eq58}), we present
\begin{equation}
\label{eq65}
 \begin{aligned}
 \psi \left( \vartheta \right)
&=\sum\limits_{m=0}^\infty {B_m \vartheta
^{2m}}
=\psi \left( 0 \right)\prod\limits_{\vartheta _k } {\left(
{1-\frac{\vartheta }{\vartheta _k }} \right)} \\
&=\psi \left( 0 \right)\prod\limits_{\Im \left( {\vartheta _k } \right)>0}
{\left( {1-\frac{\vartheta ^2}{\vartheta _k^2 }} \right)} \\
&=\psi \left( 0
\right)\prod\limits_{\Im \left( {\vartheta _k } \right)>0} {\left(
{1+\frac{\vartheta ^2}{\sigma _k^2 }} \right)} \\
&=\psi \left( 0 \right)\prod\limits_{k=1}^\infty {\left( {1+\frac{\vartheta
^2}{\sigma _k^2 }} \right)} . \\
\end{aligned}
\end{equation}

\subsection{Prove that all zeros of $\aleph \left( \tau \right)$ are real}

Taking $\vartheta =i\tau $ in (\ref{eq65}) implies that
\begin{equation}
\label{eq66}
\begin{array}{l}
 \psi \left( {i\tau } \right)
 =\sum\limits_{m=0}^\infty {B_m \left( {i\tau }
\right)^{2m}}
 =\psi \left( 0 \right)\prod\limits_{k=1}^\infty {\left[ {1+\frac{\left(
{i\tau } \right)^2}{\sigma _k^2 }} \right]} . \\
 \end{array}
\end{equation}
To simplify (\ref{eq66}), we obtain
\begin{equation}
\label{eq67}
\begin{array}{l}
 \psi \left( {i\tau } \right)
 =\sum\limits_{m=0}^\infty {\left( {-1}
\right)^mB_m \tau ^{2m}}
 =\psi \left( 0 \right)\prod\limits_{k=1}^\infty {\left( {1-\frac{\tau
^2}{\sigma _k^2 }} \right)} .
 \end{array}
\end{equation}
Taking
\begin{equation}
\label{eq68}
\aleph \left( \tau \right)=\psi \left( {i\tau } \right)
\end{equation}
in (\ref{eq67}), we have
\begin{equation}
\label{eq69}
\begin{array}{l}
 \aleph \left( \tau \right)
 =\sum\limits_{m=0}^\infty {\left( {-1}
\right)^mB_m \tau ^{2m}}
=\psi \left( 0 \right)\prod\limits_{k=1}^\infty
{\left( {1-\frac{\tau ^2}{\sigma _k^2 }} \right)}
 =\aleph \left( 0 \right)\prod\limits_{k=1}^\infty {\left( {1-\frac{\tau
^2}{\sigma _k^2 }} \right)} ,
 \end{array}
\end{equation}
which leads to
\begin{equation}
\label{eq70}
\begin{array}{l}
 \aleph \left( 0 \right)=B_0
 =\psi \left( 0 \right)
 \end{array}
\end{equation}
when one substitutes $\tau =0$ into (\ref{eq69}).

Similarly, combining (\ref{eq65}), (\ref{eq68}) and (\ref{eq70}), we have
\begin{equation}
\label{eq71}
\begin{array}{l}
\psi \left( {i\tau } \right)=\psi \left( 0 \right)\prod\limits_{\vartheta _k
} {\left( {1-\frac{\tau i}{\vartheta _k }} \right)}
\end{array}
\end{equation}
such that
\begin{equation}
\label{eq72}
\begin{array}{l}
 \aleph \left( \tau \right)
 =\psi \left( 0 \right)\prod\limits_{\vartheta _k
} {\left( {1-\frac{\tau i}{\vartheta _k }} \right)}
 =\aleph \left( 0 \right)\prod\limits_{\vartheta _k } {\left( {1-\frac{\tau
i}{\vartheta _k }} \right)} .
 \end{array}
\end{equation}
Since (\ref{eq69}) is equivalent to (\ref{eq72}), we have the identity
\begin{equation}
\label{eq73}
\begin{array}{l}
 \aleph \left( \tau \right)
 =\sum\limits_{m=0}^\infty {\left( {-1}
\right)^mB_m \tau ^{2m}}
 =\aleph \left( 0 \right)\prod\limits_{\vartheta _k } {\left( {1-\frac{\tau
i}{\vartheta _k }} \right)}
 =\aleph \left( 0 \right)\prod\limits_{k=1}^\infty {\left( {1-\frac{\tau
^2}{\sigma _k^2 }} \right)} .
 \end{array}
\end{equation}
From (\ref{eq73}) it is observed that all zeros of $\aleph \left( \tau \right)$
are real.

Thus, this is required result.

\begin{remark}
Replacing $\overline \vartheta \in \mathbb{C} $ by $\vartheta \in \mathbb{C}$ in (\ref{eq35}), we also deduce that
\[
\begin{array}{l}
\psi \left( {\vartheta
} \right)
 =\psi \left( 0 \right)\prod\limits_{\Im \left( {\vartheta _k } \right)>0}
{\left( {1-\frac{\vartheta ^2}{\vartheta _k^2 }} \right)}
 =\psi \left( 0 \right)\prod\limits_{\Im \left( {\vartheta _k } \right)>0}
{\left( {1-\frac{\vartheta ^2}{\vartheta _k^2 }}
\right)},
 \end{array}
\]
which leads to (\ref{eq43}) and (\ref{eq48}) when $\vartheta=t \in \mathbb{R}$.

\end{remark}

\subsection{Equivalently sufficient conditions}

With use of Lemma \ref{Le1} and (\ref{eq73}), we obtain
\begin{equation}
\label{eq74}
 \begin{aligned}
 \aleph \left( \tau \right)
&=\sum\limits_{m=0}^\infty {\alpha _m \tau ^{2m}}
=\aleph \left( 0 \right)\prod\limits_{\Im \left(
{\vartheta _k } \right)>0} {\left( {1+\frac{\tau ^2}{\vartheta _k^2 }}
\right)}
 =\aleph \left( 0 \right)\prod\limits_{k=1}^\infty {\left( {1-\frac{\tau
^2}{\sigma _k^2 }} \right)} \\
&=\aleph \left( 0 \right)\prod\limits_{\vartheta
_k } {\left( {1-\frac{\tau i}{\vartheta _k }} \right)}.
 \end{aligned}
\end{equation}
From (\ref{eq74}) we know that $\pm \sigma _k \left( {\sigma _k >0} \right)$ are all real zeros of
$\aleph \left( \tau \right)$.

As a direct result of (\ref{eq74}), we have the following:

\begin{corollary}
\label{CO1}
If $\aleph \in \wp $, then there exist the following equivalent
representations:

(A) All of the zeros of $\aleph \left( \tau \right)$ are real.

(B) There exists the identity
\begin{equation}
\label{eq75}
\begin{array}{l}
\sum\limits_{m=0}^\infty {\alpha _m \tau ^{2m}} =\aleph
\left( 0 \right)\prod\limits_{\Im \left( {\vartheta _k } \right)>0} {\left(
{1+\frac{\tau ^2}{\vartheta _k^2 }} \right)}.
\end{array}
\end{equation}

(C) There exists the identity
\begin{equation}
\label{eq76}
\begin{array}{l}
\sum\limits_{m=0}^\infty {\alpha _m \tau ^{2m}} =\Lambda
\left( 0 \right)\prod\limits_{\vartheta _k } {\left( {1-\frac{\tau
i}{\vartheta _k }} \right)}.
\end{array}
\end{equation}

(D) There exists the identity
\begin{equation}
\label{eq77}
\begin{array}{l}
\sum\limits_{m=0}^\infty {\alpha _m \tau ^{2m}} =\Lambda
\left( 0 \right)\prod\limits_{k=1}^\infty {\left( {1-\frac{\tau ^2}{\sigma
_k^2 }} \right)}.
\end{array}
\end{equation}
\end{corollary}

By Corollary \ref{CO1}, we see that there are some sufficient conditions that all
of the zeros of $\aleph \left( \tau \right)$ are real.

\section{A typical application}{\label{sec:4}}

\subsection{Prove the conjecture of Lagarias and Montague}
Let us start with its proof.
By using (\ref{eq5}), the Lagarias-Montague function  $\Phi \left( \tau\right)$ can be written as
\begin{equation}
\label{eq78}
\begin{array}{l}
 \Phi \left( \tau \right)
 =2\int\limits_0^\infty {g\left( \ell \right)\ell
^{-1}\sin \left( {\tau\ell } \right)d\ell }
 =2\int\limits_0^\infty {g\left( \ell \right)\ell ^{-1}\left[
{\sum\limits_{m=0}^\infty {\frac{\left( {-1} \right)^n\left( {\tau\ell }
\right)^{2m+1}}{\left( {2m+1} \right)!}} } \right]d\ell } .
 \end{array}
\end{equation}
To simplify (\ref{eq78}), we obtain
\begin{equation}
\label{eq79}
\Phi \left( \tau \right)=\sum\limits_{m=0}^\infty {\left( {-1}
\right)^n\mathchar'26\mkern-10mu\lambda _{2m+1} \tau^{2m+1}} ,
\end{equation}
where
\begin{equation}
\label{eq80}
\begin{array}{l}
\mathchar'26\mkern-10mu\lambda _{2m+1} =2\int\limits_0^\infty {g\left( \ell
\right)\frac{\ell ^{2m}}{\left( {2n+1} \right)!}d\ell } >0.
\end{array}
\end{equation}
Because of
\begin{equation}
\label{eq81}
\begin{array}{l}
\Xi \left( \tau \right)=\Phi ^{\left( 1 \right)}\left(\tau \right),
\end{array}
\end{equation}
In fact, Boas (see Theorem 2.4.1 in \cite{2}, p.13)  argued that $\Phi \left(
\tau \right)$ and $\Xi \left( \tau \right)$ are of the same order and type. In
view of the work of Dimitrov and Lucas \cite{20}, $\Phi \left( \tau \right)$ and $\Xi
\left(\tau \right)$ are the functions order $\rho _1 =1$.

Let
\begin{equation}
\label{eq82}
\begin{array}{l}
\widehat{\Phi }_1 \left( \tau \right)=\frac{\Phi \left( {i\tau}
\right)}{i\tau}=\sum\limits_{m=0}^\infty {\mathchar'26\mkern-10mu\lambda
_{2m+1} \tau^{2m}}
\end{array}
\end{equation}
and
\begin{equation}
\label{eq83}
\begin{array}{l}
\widehat{\Phi }_2 \left( \tau \right)=\frac{\Phi \left(\tau
\right)}{\tau}=\sum\limits_{m=0}^\infty {\left( {-1}
\right)^m\mathchar'26\mkern-10mu\lambda _{2m+1} \tau^{2m}}
\end{array}
\end{equation}
such that
\begin{equation}
\label{eq84}
\begin{array}{l}
F\left( \tau \right)=i\tau\widehat{\Phi }_1 \left(\tau
\right)
\end{array}
\end{equation}
and
\begin{equation}
\label{eq85}
\begin{array}{l}
\Phi \left( \tau \right)=\tau \widehat{\Phi }_2 \left(\tau \right).
\end{array}
\end{equation}
Since (\ref{eq82}) and (\ref{eq83}) are of the same order and type due to the fact
\[
\mathchar'26\mkern-10mu\lambda _{2m+1} =\left| {\left( {-1}
\right)^m\mathchar'26\mkern-10mu\lambda _{2m+1} } \right|,
\]
(\ref{eq84}) and (\ref{eq85})
are also of the same order and type. This implies that (\ref{eq82}) , (\ref{eq83}), (\ref{eq84})
and (\ref{eq85}) are of order $\rho _1 =1$. Moreover, (\ref{eq82}) and (\ref{eq83}) are the even entire
functions.

Since (\ref{eq82}) is an even entire function of order $\rho _1 =1$ with the
positive real coefficients $\mathchar'26\mkern-10mu\lambda _{2m+1} >0$,
Theorem 3 in Levin's book (see \cite{1}, p.8) said that the product presentation of
$\Phi \left( {it} \right)$ reads
\begin{equation}
\label{eq86}
\begin{array}{l}
F\left( \tau \right)=it\widehat{\Phi }_1 \left(\tau \right)=i\tau\widehat{\Phi }_1
\left( 0 \right)\prod\limits_{u_k } {\left( {1-\frac{\tau}{u_k }} \right)}
\end{array}
\end{equation}
with
\begin{equation}
\label{eq87}
\begin{array}{l}
\widehat{\Phi }_1 \left( \tau \right)=\widehat{\Phi }_1 \left( 0
\right)\prod\limits_{u_k } {\left( {1-\frac{\tau}{u_k }} \right)} ,
\end{array}
\end{equation}
where $u_k $ run the zeros of $\widehat{\Phi }_1 \left(\tau \right)$.

If $\Im \left( {u_k } \right)>0$ and (\ref{eq82}) is an even function with the
complex zeros, then there is the functional equation
\begin{equation}
\label{eq88}
\widehat{\Phi }_1 \left( \tau \right)=\widehat{\Phi }_1 \left( {-\tau} \right)
\end{equation}
such that
\begin{equation}
\label{eq89}
\begin{array}{l}
 \widehat{\Phi }_1 \left( \tau \right)
 =\widehat{\Phi }_1 \left( 0
\right)\prod\limits_{u_k } {\left( {1-\frac{\tau}{u_k }} \right)}\\
=\widehat{\Phi }_1 \left( 0 \right)\prod\limits_{\Im \left( {u_k }
\right)>0} {\left( {1-\frac{\tau}{u_k }} \right)\left( {1+\frac{\tau}{u_k }}
\right)} \\
 =\widehat{\Phi }_1 \left( 0 \right)\prod\limits_{\Im \left( {u_k }
\right)>0} {\left( {1-\frac{\tau^2}{u_k^2 }} \right)} . \\
 \end{array}
\end{equation}
From (\ref{eq82}) and (\ref{eq83}) we have
\begin{equation}
\label{eq90}
\widehat{\Phi }_1 \left( {it} \right)=\widehat{\Phi }_2 \left( t \right)
\end{equation}
and
\begin{equation}
\label{eq91}
\begin{array}{l}
\widehat{\Phi }_1 \left( 0 \right)=\widehat{\Phi }_2 \left( 0
\right)=2\int\limits_0^\infty {g\left( \ell \right)\frac{\ell ^{2m}}{\left(
{2n+1} \right)!}d\ell } >0.
\end{array}
\end{equation}
With the aid of (\ref{eq89}), (\ref{eq90}) and (\ref{eq91}), we may get
\begin{equation}
\label{eq92}
\begin{array}{l}
 \widehat{\Phi }_2 \left( \tau \right)
 =\widehat{\Phi }_1 \left( 0
\right)\prod\limits_{\Im \left( {u_k } \right)>0} {\left(
{1+\frac{\tau^2}{u_k^2 }} \right)}
 =\widehat{\Phi }_2 \left( 0 \right)\prod\limits_{\Im \left( {u_k }
\right)>0} {\left( {1+\frac{\tau^2}{u_k^2 }} \right)} . \\
 \end{array}
\end{equation}
By combination of (\ref{eq83}) and (\ref{eq92}), we present
\begin{equation}
\label{eq93}
\begin{array}{l}
 \widehat{\Phi }_2 \left(\tau \right)=\sum\limits_{m=0}^\infty {\left( {-1}
\right)^m\mathchar'26\mkern-10mu\lambda _{2m+1} \tau^{2m}}
 =\widehat{\Phi }_2 \left( 0 \right)\prod\limits_{\Im \left( {u_k }
\right)>0} {\left( {1+\frac{\tau^2}{u_k^2 }} \right)} . \\
 \end{array}
\end{equation}
Because of the fact $\widehat{\Phi }_2 \left( \tau \right)$ is of order $\rho _1
=1$, Theorem 13 (\cite{1}, p.24) has reported that there exists  any $\widehat{\varepsilon }>0$ such that the series
\begin{equation}
\label{eq94}
\sum\limits_{k=1}^\infty {\frac{1}{|u_k|^{1+\widehat{\varepsilon }} }}
\end{equation}
is convergent.

Taking $\widehat{\varepsilon }=1$ in (\ref{eq94}), the series
\begin{equation}
\label{eq95}
\sum\limits_{k=1}^\infty {\frac{1}{|u_k|^2 }}
\end{equation}
is convergent and (\ref{eq92}) is also convergent.

In sum, we have the followings three conditions:

(A1) $\widehat{\Phi }_2 \left( \tau \right)$ is of order $\rho _1 =1$.

(A2) The identity (\ref{eq93}) holds for $\tau\in \mathbb{C}$.

(A3)The series (\ref{eq95}) is convergent.

Then, we obtain
\begin{equation}
\label{eq96}
\widehat{\Phi }_2 \in \wp .
\end{equation}
By using Theorem \ref{Th1} and (B) in Corollary \ref{CO1}, we have
from (\ref{eq93}) that all zeros $u_k $ of the even entire function $\widehat{\Phi
}_2 \left( \tau \right)$ are real if $\widehat{\Phi }_2 \in \wp $.

This implies that
\begin{equation}
\label{eq97}
u_k =\pm i\beta _k ,
\end{equation}
where $\beta _k >0$.

By using (\ref{eq97}), the identity (\ref{eq93}) can be written as
\begin{equation}
\label{eq98}
\begin{array}{l}
 \widehat{\Phi }_2 \left( \tau \right)=\sum\limits_{m=0}^\infty {\left( {-1}
\right)^m\mathchar'26\mkern-10mu\lambda _{2m+1} \tau^{2m}}
 =\widehat{\Phi }_2 \left( 0 \right)\prod\limits_{\Im \left( {u_k }
\right)>0} {\left( {1-\frac{\tau^2}{\beta _k^2 }} \right)} . \\
 \end{array}
\end{equation}
With (\ref{eq79}), (\ref{eq83}) and (\ref{eq98}), we clearly see that
\begin{equation}
\label{eq99}
 \begin{aligned}
 \Phi \left( \tau \right)
 &=\sum\limits_{m=0}^\infty {\left( {-1}
\right)^n\mathchar'26\mkern-10mu\lambda _{2m+1} \tau^{2m+1}} \\
&=\tau\widehat{\Phi
}_2 \left( 0 \right)\prod\limits_{\Im \left( {u_k } \right)>0} {\left(
{1-\frac{\tau^2}{\beta _k^2 }} \right)} \\
 &=\tau\widehat{\Phi }_2 \left( 0 \right)\prod\limits_{\beta _k >0} {\left(
{1-\frac{\tau^2}{\beta _k^2 }} \right)} \\
&=\tau\widehat{\Phi }_2 \left( 0
\right)\prod\limits_{k=1}^\infty {\left( {1-\frac{\tau^2}{\beta _k^2 }}
\right)} . \\
 \end{aligned}
\end{equation}
It follows from (\ref{eq99}) that all zeros of
$\Phi \left(\tau \right)$ are real because $\beta _k >0$.

Thus, we prove the conjecture of Lagarias and Montague.

\begin{remark}
By substitution of (\ref{eq97}) into (\ref{eq89}), we have
\begin{equation}
\label{eq100}
 \begin{aligned}
 \widehat{\Phi }_1 \left(\tau \right) &=\sum\limits_{m=0}^\infty
{\mathchar'26\mkern-10mu\lambda _{2m+1} \tau^{2m}} \\
&=\widehat{\Phi }_1 \left( 0
\right)\prod\limits_{u_k } {\left( {1-\frac{\tau}{i\beta _k }} \right)} \\
 &=\widehat{\Phi }_1 \left( 0 \right)\prod\limits_{\Im \left( {u_k }
\right)>0} {\left( {1-\frac{\tau}{i\beta _k }} \right)\left( {1+\frac{\tau}{i\beta
_k }} \right)} \\
&=\widehat{\Phi }_1 \left( 0
\right)\prod\limits_{k=1}^\infty {\left( {1+\frac{\tau^2}{\beta _k^2 }}
\right)} , \\
 \end{aligned}
\end{equation}
which implies that
\begin{equation}
\label{eq101}
\begin{array}{l}
 F\left(\tau \right)
=i\tau\sum\limits_{m=0}^\infty {\mathchar'26\mkern-10mu\lambda _{2m+1}
\tau^{2m}}
 =i\tau\widehat{\Phi }_1 \left( 0 \right)\prod\limits_{u_k } {\left(
{1-\frac{\tau}{u_k }} \right)}
=i\tau\widehat{\Phi }_1 \left( 0
\right)\prod\limits_{k=1}^\infty {\left( {1+\frac{\tau^2}{\beta _k^2 }}
\right)} .
 \end{array}
\end{equation}
\end{remark}

\subsection{Two open problems}
From (\ref{eq100}) and (\ref{eq101}) it is easily seen that all zeros of $\widehat{\Phi }_1 \left( \tau
\right)$ are $\tau=\pm i\beta _k $, where $\beta _k >0$, and that the function
$F\left( \tau \right)$ has the purely imaginary number zeros $\tau=\pm i\beta _k
$, where $\beta _k >0$, and real zero $t=0$.
Similarly, by (\ref{eq98}) and (\ref{eq99}), it is also observed that all real zeros of $\widehat{\Phi }_2
\left( \tau \right)$ are $\tau=\pm \beta _k $, where $\beta _k >0$, and that the
function $\Phi \left( \tau \right)$ has the real zeros $\tau=\pm \beta _k $, where
$\beta _k >0$ and real zero $\tau=0$. Here, we call $\widehat{\Phi }_1 \left( \tau
\right)$ as the hungry pair of $\widehat{\Phi }_2 \left( \tau
\right)$ if there exist all zeros $\tau=\pm \beta _k $ of $\widehat{\Phi }_2 \left( \tau
\right)$ and all zeros $\tau=\pm i\beta _k $ of $\widehat{\Phi }_1 \left( \tau
\right)$. Computing real zeros of $(\ref{eq5})$ is still an open problem in the theory of the Lagarias-Montague
function.

By using the observation that all zeros of
\[
\Phi \left(t
\right)=\int\limits_0^t {\Xi \left( t \right)dt}
\]
 and $\Xi \left( t\right)=\Phi ^{\left( 1 \right)}\left(t \right)$
 are real, we have the followings:

\begin{problem}
\label{PR1}
Let $\aleph \in \wp $. Then all zeros of $G\left(\tau
\right)$ are real if $G\left(t
\right)=\aleph ^{\left( 1 \right)}\left(t\right)$.
\end{problem}
As an analogous problem \ref{PR1}, the real zeros of the derivative of
the entire function in the Laguerre-P\'{o}lya class was proposed by P\'{o}lya in 1913 \cite{19} and
proved by Hellerstein and Williamson \cite{20,21}.
\begin{problem}
\label{PR2}
Let $\aleph \in \wp $. Then all zeros of $\mathbb{ M}\left(\tau \right)$ are real if there exists the integral
\begin{equation}
\label{eq102}
\mathbb{ M}\left( t \right)=\int\limits_0^t {\aleph \left( t \right)dt}.
\end{equation}
\end{problem}

Here, we easily find that all
zeros of the function $\mathbb{ M}\left(\tau \right)=\cos \left( \tau \right)$ are real if
\[
\cos \left( t \right)=\int\limits_0^t {\sin \left( t \right)dt},
\]
and that all zeros of the function $G\left(\tau
\right)=-\sin \left( \tau \right)$ are
real if $\cos ^{\left( 1 \right)}\left( t \right)=-\sin \left( t \right)$ .
Also, it is easy to see that $cosh \left( \tau
\right)$ is considered as the hungry pair of $cos \left( \tau
\right)$.
As a direct result of Corollary \ref{CO1}, we have the following:
\begin{corollary}
\label{CO2}
$cos \left( \tau\right)$ belongs to the class $\wp$.
\end{corollary}
\begin{proof}
Adopting the product representations of $cos \left( \tau\right)$ (see \cite{22}, p.114), we structure
the class of the series and product representations of $cos \left( \tau\right)$, given as
\begin{equation}
\label{eq103}
\begin{array}{l}
\sum\limits_{m = 0}^\infty  {\frac{{\left( { - 1} \right)^m }}{{\left( {2m} \right)!}}\tau ^{2m} }
= \prod\limits_{k = 1}^\infty  {\left[ {1 - \frac{{\tau ^2 }}{{\left( {k - \frac{1}{2}} \right)^2 \pi ^2 }}} \right]}.
\end{array}
\end{equation}
By Corollary \ref{CO1} and (\ref{eq103}), we directly obtain the required result because
$cos \left( \tau\right)$ is an even function of order $\gamma =1$ (see \cite{22}, p.255) and the Euler's product of
(\ref{eq103}) is convergent.
\end{proof}
\section{Conclusion} {\label{sec:5}}

In the present article we have proposed a sufficient condition for all real
zeros of a class of the even entire function. With the aid of the obtained result, we have
proved that the conjecture of Lagarias and Montague is true. By comparison
between the real zeros of the Riemann Xi  and the Lagarias-Montague
functions, we have suggested two open problems for the even entire functions.
The result may be proposed as a new mathematical approach to open a new door
for handling the de Bruijn-Newman constant.

\end{document}